\documentclass[11pt]{article}

\usepackage{fullpage}
\usepackage{standard}
\usepackage{groups}
\usepackage{natbib}

\usepackage{xypic}	
\usepackage[all,2cell,ps]{xy}
\usepackage{epsfig}
\usepackage{young}
\usepackage[enableskew]{youngtab}

\newtheorem{thm}{Theorem}

\newenvironment{pf}{\textbf{Proof.}}{\m{\blacksquare}}
\renewcommand{\marginpar}[1]{}

\newcommand{\Lslice}[1]{\!\nts \downharpoonleft_{#1}\!\!}
\newcommand{\Rslice}[1]{\!\nts \downharpoonright_{#1}\!\!}
\newcommand{\absG}{\vert\tts G\tts \vert}
\newcommand{\downarr}{\!\!\downarrow}
\newcommand{\mirr}[1]{{#1}^{-}}

\begin{document}

\title{The skew spectrum of functions on finite groups and their homogeneous spaces}
\author{
{\bf Risi Kondor}\\
\texttt{risi@gatsby.ucl.ac.uk}\\
Gatsby Unit for Computational Neuroscience\\
University College London\\
17 Queen Square, WC1N 3AR\\
United Kingdom
}
\date{}
\maketitle
\begin{abstract}
  Whenever we have a group acting on a class of functions by translation, the 
  bispectrum offers a principled and lossless way of representing such functions invariant to the action.   
  Unfortunately, computing the bispectrum is often costly and complicated. 
  In this paper we propose a unitarily equivalent, but 
  easier to compute set of invariants, which we call the skew spectrum. 
  For  functions on homogeneous spaces 
  the skew spectrum can be efficiently computed using some ideas from Clausen-type 
  fast Fourier transforms.     
\end{abstract}

\section{Introduction}

Given a function \m{f\colon \RR \to \CC} with Fourier transform \m{\h f}, 
it is well known that the bispectrum 
\begin{equation}\label{eq: bispectrum Euclidean}
  \h q(k_1,k_2)=\h f(k_1)\, \h f(k_2)\,\h f(k_1+k_2) \qqquad k_1,k_2\in \RR
\end{equation}
is invariant to translations of \m{f}. 
The bispectrum is used in various signal and image processing applications which 
involve translation invariance. 

Kakarala generalized \rf{eq: bispectrum Euclidean} to signals on a wide class of non-commutative groups, 
including all compact, and thus, all finite groups. 
On a group \m{G} 
the (left-) translate of \m{f\colon G\to\CC} is defined \m{f^t(x)=f(t^{-1}x)}, and 
the bispectrum takes the form 
\begin{equation}\label{eq: bispectrum1}
  \h q(\rho_1,\rho_2)=\brbig{\h f(\rho_1) \otimes \h f(\rho_2)}^\dag\: 
  C_{\rho_1\rho_2}
  \sqbbigg{\:\bigoplus_{\rho\in\Rcal} {\h f(\rho)}^{\oplus m_\rho} \:}  C_{\rho_1\rho_2}^\dag, 
\end{equation}
where \m{\rho_1,\rho_2} and \m{\rho} are irreducible representations of \m{G}, 
and \m{C_{\rho_1\rho_2}} are the Clebsch-Gordan matrices. 
Recall that the Fourier transform on a non-commutative group 
is a collection of matrices 
indexed by the irreducible representations.   
The remarkable fact shown by Kakarala is that for a wide class of groups 
the bispectrum is not only invariant to left translation, but,  
assuming that each Fourier component \m{\h f(\rho)} is invertible, 
the \m{\h q(\rho_1,\rho_2)} matrices also uniquely determine \m{f}  
up to translation \citep[Theorem 3.2.4]{KakaralaPhD}. 
In other words, given \m{f,f'\colon G\to\CC} and their respective bispectra 
\m{\br{\h q(\rho_1,\rho_2)}_{\rho_1,\rho_2}} and  \m{\br{\h q'(\rho_1,\rho_2)}_{\rho_1,\rho_2}}, we have  
\m{q(\rho_1,\rho_2)=q'(\rho_1,\rho_2)} if and only if \m{f'(x)=f(t^{-1}x)} for some \m{t\in G}. 

Unfortunately, the bispectrum can be difficult to compute, even on finite groups of moderate size. 
The total number of entries in the \m{\h q(\rho_1,\rho_2)} matrices is \m{\abs{G}^2}, which can be a problem. 
Another issue is the difficulty of deriving the form of the Clebsch Gordan matrices for non-trivial groups. 
Finally, computing \rf{eq: bispectrum1} is expensive because it involves the multiplication of large matrices. 
These factors all limit the practical applicability of the bispectrum method. 

In this paper we introduce the \df{skew spectrum}, which furnishes 
an alternative set of invariants.
The skew spectrum can be computed by Fourier transformation alone, 
without any explicit Clebsch-Gordan transforms. 
We show that the skew spectrum is unitarily equivalent to the bispectrum. 
In particular, the skew spectrum is sufficient to uniquely determine the original function \m{f} up to translation.  

An important generalization of the above is when \m{f} is not a function on \m{G} itself, but 
a function on a homogeneous space of \m{G}. 
This case is likely to be more relevant to applications than the original case of functions on \m{G}.  
The new invariants can be restricted to a wide range of homogeneous spaces,  
and are closely adapted to Clausen-type FFTs, offering computational economy. 


\ignore{
Given a finite group \m{G} and a function \m{f\colon G\to \CC}, the \df{left-translate} of \m{f} by \m{g\<\in G} is 
the function \m{f^g(x)=f(g^{-1} x)}. 
More generally, if \m{G} acts transitively on a set \m{X}, we define the left-translate of \m{f\colon X\to \CC} as 
\m{f^g(x)=f(g^{-1}(x))}. 

A functional \m{p\colon f\mapsto p(f)} is said to be \df{invariant} to left translation 
if \m{p(f)=p(f^g)} for any \m{f\in L_2(X)} and any \m{g\in G}. 
A set \m{P=\cbr{\seq{p}{k}}} of such invariants is said to be \df{complete} if 
\m{p_1(f)\<=p_1(f'),\:\ldots,\: p_k(f)\<=p_k(f')} implies \m{f=f^g} for some \m{g\in G}. 
A complete set \m{P} of invariant functionals is said to be \df{minimal}, if omitting any single functional from \m{P}, 
the remaining set is no longer complete.

In this paper we construct minimal complete sets of invariants. 
Each \m{p\in P} will be a third order polynomial in the variables \m{\cbr{f(x)}_{x\in X}}. 
We show that \m{P} is complete, but in general not minimal. 
}

\section{Fourier transforms and the bispectrum}\label{sec: bispectrum}

In this paper \m{G} denotes a finite group of cardinality \m{\abs{G}}  
and \m{\Rcal} denotes a complete set of inequivalent 
irreducible complex-valued matrix representations of \m{G}. 
Since \m{G} is finite, \m{\Rcal} is a finite set. 
We use the letter \m{\rho} to denote the individual irreducibles, and \m{d_\rho} to denote their dimensionalities. 
Thus, each \m{\rho\in\Rcal} is a homomorphism \m{\rho\colon G\to \CC^{d_\rho\times d_\rho}}.   
Without loss of generality, 
we assume that all \m{\rho\in\Rcal} are unitary, i.e., \m{\rho(x^{-1})=\br{\rho(x)}^{-1}=\rho(x)^\dag}. 

We \marginpar{translation} denote the set of complex valued functions on \m{G} by \m{L(G)}. 
Given \m{f\in L(G)} and a group element \m{t}, the \df{left-translate} of \m{f} by \m{t} is the function 
\m{f^t\in L(G)} defined \m{f^t(x)=f(t^{-1}x)}. 
A functional \m{p\colon f\mapsto p(f)} is said to be \df{invariant} to left translation 
(or just \df{left-invariant} for short) 
if \m{p(f)=p(f^t)} for any \m{f\in L(G)} and any \m{t\in G}. 
A set \m{P=\cbr{\seq{p}{k}}} of such invariants is said to be \df{complete} if 
\m{p_1(f)\<=p_1(f'),\:\ldots,\: p_k(f)\<=p_k(f')} implies \m{f'=f^t} for some \m{t\in G}. 
The objective of this paper is to find an efficiently computable yet complete set of left-invariant functionals. 

\subsection{The Fourier transform}

The \marginpar{Fourier transform} \df{Fourier transform} of a function \m{f\in L(G)} is the collection of matrices 
\begin{equation*}
  \h f(\rho)=\sum_{x \in G} \rho(x)\, f(x) \qqquad\qqquad \rho\in\Rcal.  
\end{equation*}
With respect to the norms 
\begin{equation*}
  \nm{f}^2=\ovr{\abs{G}} \sum_{x\in G}\abs{f(x)}^2 
  \qquad\text{and}\qquad
  \nmbig{\h f}^2=\ovr{\abs{G}^2}\sum_{\rho\in\Rcal} d_\rho \nmbig{\h f(\rho)}_{2}^2, 
\end{equation*}
(where \m{\nm{M}_2} denotes the Frobenius norm  
\m{\brbig{\sum_{i,j=1}^n \sqb{M}_{i,j}}^{1/2}}~)
the transform \m{\mathfrak{F}\colon f\mapsto\h f} is unitary, and its inverse is given by 
\begin{equation*}
  f(x)=\ovr{\abs{G}} \sum_{\rho\in\Rcal} d_\rho\: \tr\brbig{\rho(x) \cdot \h f(\rho)} \qqquad \qqquad x\in G. 
\end{equation*}
Defining the \df{convolution} of two functions \m{f,g\in L(G)} as 
\begin{equation*}
  \br{f\ast g}(x)=\sum_{y\in G} f(xy^{-1})\: g(y),
\end{equation*}
it is easy to verify 
the convolution theorem \m{\widehat{f \ast g}(\rho)=\h f(\rho) \cdot \h g(\rho)}. 
Another property of interest   
is the behavior of the Fourier transform under translation: 
\begin{equation}\label{eq: translation}
  \h {f^t}(\rho)=\sum_{x \in G} \rho(x)\, f(t^{-1}x)
  =\sum_{x \in G} \rho(tx)\, f(x)
  =\rho(t)\sum_{x \in G} \rho(x)\, f(x)
  =\rho(t)\cdot \h f(\rho). 
\end{equation}
Finally, defining \m{\mirr{f}\in L(G)} by \m{{\mirr{f}}(x)=f(x^{-1})}, 
and letting \m{{}^\ast} denote complex conjugation,  
we have \m{\h {{\mirr{f}}}(\rho)^\ast=\h f(\rho)^\dag}.  

In \marginpar{covariance}
general, we say that a matrix-valued functional \m{F\colon L(G) \to \CC^{D \times d_\rho}} is \m{\rho}-\df{covariant} 
if \m{F(f^t)=\rho(t)\cdot  F(f)} for any \m{f\in L(X)} and any \m{t\in G}. 
Similarly, \m{F'\colon L(X) \to \CC^{d_\rho\times d_\rho}} is \m{\rho}-\df{contravariant} if 
\m{F'(f^t)=F'(f)\cdot \rho(t)^\dag}. 
Clearly, the adjoint (conjugate transpose) of a \m{\rho}-covariant functional is \m{\rho}-contravariant. 

\subsection{The power spectrum}

If \marginpar{spectrum} 
\m{F} is \m{\rho}-covariant and \m{F'} is \m{\rho}-contravariant, then 
\m{Q(f)=F'(f)\: F(f)} is invariant, i.e., \m{Q(f^t)=Q(f)}. 
Thus, 
the matrices \m{\h a(\rho)=\h f(\rho)^\dagger \: \h f(\rho)} are a natural starting point for finding 
left-invariant functionals. 
On closer inspection, the invariance of \m{\h a} should not come as a surprise, since 
by the convolution theorem \m{\h a}  
is the Fourier transform of the \df{autocorrelation} 
\begin{equation*}
  a(x)=
        \sum_{y\in G} \mirr{f}(xy^{-1})^\ast\,f(y)=
        \sum_{y\in G} f(yx^{-1})^\ast\,f(y),
\end{equation*}
which is manifestly invariant to left-translation. 
We call \m{\br{\h a(\rho)}_{\rho\in\Rcal}} the \df{power spectrum} of \m{f}. 

Unfortunately, the matrix elements of \m{\br{\h a(\rho)}_{\rho\in\Rcal}} 
fall short of forming a complete set of invariants.  
This is because the spectrum is insensitive to the relative 
``phase'' of the various Fourier components.  
The natural way to couple the components is to take tensor products, 
forming higher order spectra. 

\subsection{Triple correlation and the bispectrum}

Recall \marginpar{Clebsch-Gordan}
that for any \m{\rho_1, \rho_2\in\Rcal}, the tensor product representation 
\m{\rho_1\otimes \rho_2} decomposes into irreducible components in the form 
\begin{equation}\label{eq: deco1}
  (\rho_1 \otimes \rho_2)(x)=
  \rho_1(x)\otimes \rho_2(x)=
  C_{\rho_1\rho_2}\,
  \sqbbigg{\:\bigoplus_{\rho\in\Rcal} \rho(x)^{\oplus m_\rho}\:}\:  C_{\rho_1\rho_2}^\dag, 
\end{equation}
where \m{C_{\rho_1\rho_2}} is a unitary matrix called the \df{Clebsch-Gordan} matrix, 
\m{m_\rho=m(\rho_1,\rho_2,\rho)\in \ZZ^{+}} is the multiplicity of \m{\rho} in the decomposition, and 
\m{M^{\oplus k}} is a shorthand for \m{\bigoplus_{i=1}^k M}.  

Now \marginpar{bispectrum} 
consider the tensor product of the Fourier matrices \m{\h f(\rho_1)} and \m{\h f(\rho_2)}. 
Under translation 
\begin{equation*}
  \h f(\rho_1)\otimes \h f(\rho_2)\quad \mapsto\quad  
  \h {f^t}(\rho_1)\otimes \h {f^t}(\rho_2)=
  \br{\rho_1(t)\otimes \rho_2(t)} \cdot \brN{\h f(\rho_1)\otimes \h f(\rho_2)}, 
\end{equation*}
hence
\begin{equation*}
  C_{\rho_1\rho_2}^\dag\ts \brbig{\h {f^t}(\rho_1)\otimes \h {f^t}(\rho_2)}\, C_{\rho_1\rho_2} =
 \sqbbigg{\:\bigoplus_{\rho\in\Rcal} \rho(t)^{\oplus m_\rho}\:}~
  C_{\rho_1\rho_2}^\dag\ts \brbig{\h f(\rho_1)\otimes \h f(\rho_2)}\, C_{\rho_1\rho_2}.    
\end{equation*}
We say that \m{C_{\rho_1\rho_2}^\dag\ts \brbig{\h {f}(\rho_1)\otimes \h {f}(\rho_2)}\, C_{\rho_1\rho_2}} 
is \m{\bigoplus_\rho \rho^{\oplus m_\rho}}-covariant. 
In particular, 
the \df{bispectrum}  
\begin{equation}\label{eq: bispectrum}
  \h b(\rho_1,\rho_2)=\brbig{\h f(\rho_1) \otimes \h f(\rho_2)}^\dag\:
  C_{\rho_1\rho_2}\:
  \sqbbigg{\:\bigoplus_{\rho\in\Rcal} {\h f(\rho)}^{\oplus m_\rho} \:}\:  C_{\rho_1\rho_2}^\dag \qqquad \rho_1,\rho_2\in\Rcal 
\end{equation}
is invariant to left-translations of \m{f}. 
The remarkable fact proved by Kakarala for compact groups in general is that 
provided that all \m{\h f(\rho)} are invertible, the matrix elements of 
\m{\brN{\h b(\rho_1,\rho_2)}_{\rho_1,\rho_2\in\Rcal}} are not only invariant, 
but also complete \citep[Theorem 3.2.4]{KakaralaPhD}\citep{KakaralaTriple}.
\footnote{Kakarala defines the bispectrum in a slightly different form, equivalent to \rf{eq: bispectrum} up to 
unitary transforms and complex conjugation.}  
It is possible to construct even higher order spectra, but there is no need: 
the bispectrum is already complete. 

It is natural to ask what the analog of the autocorrelation is, i.e., what invariant function the bispectrum is 
the Fourier transform of. 
To answer this question we consider functions on the direct product group \m{G\times G}. 
Recall from general representation theory that the irreducible representations of 
\m{G\times G} are the tensor products of the irreducible representations of \m{G}. 
Thus, 
\begin{equation*}
  \h f(\rho_1)\otimes\h f(\rho_2) =
        \sum_{x_1\in G} \sum_{x_2\in G} \br{\rho_1(x_1)\otimes \rho_2(x_2)} f(x_1)\ts f(x_2)
\end{equation*}
in \rf{eq: bispectrum} is the Fourier transform of \m{u(x_1,x_2)=f(x_1)\,f(x_2)}, 
while 
\begin{equation*}
    C_{\rho_1\rho_2}
  \sqbbigg{\:\bigoplus_{\rho\in\Rcal} {\h f(\rho)}^{\oplus m_\rho} \:}  C_{\rho_1\rho_2}^\dag = 
  \sum_{x\in G}   C_{\rho_1\rho_2}
  \sqbbigg{\:\bigoplus_{\rho\in\Rcal} {\h \rho(x)}^{\oplus m_\rho} \:}  C_{\rho_1\rho_2}^\dag\; f(x)=
  \sum_{x\in G} \br{\rho_1(x)\otimes \rho_2(x)} f(x)
\end{equation*}
is the Fourier transform of 
\begin{equation*}
  v(x_1,x_2)=
  \begin{cases}
    ~~f(x_1)&\text{if } x_1=x_2\\
    ~~0&\text{otherwise}\;. 
  \end{cases}
\end{equation*}
Thus, by the convolution theorem (on \m{G\< \times G}), \rf{eq: bispectrum} is the Fourier transform of 
\m{u^{-}\ast v}, where \m{u^{-}(x_1,x_2)=u(x_1^{-1},x_2^{-1})}. 
In long hand, the bispectrum is the Fourier transform of  
\begin{multline}\label{eq: triple correlation}
  \sum_{y_1\in G}\sum_{y_2\in G} u^{-}(x_1 y_1^{-1}\!, x_2 y_2^{-1})^\ast\, v(y_1,y_2)=
  \sum_{y_1\in G}\sum_{y_2\in G} u(y_1 x_1^{-1}\!,y_2 x_2^{-1})^\ast\, v(y_1,y_2)\\
  =\sum_{y\in G} f(y x_1^{-1})^\ast\,f(y x_2^{-1})^\ast\,f(y)
  =b\br{x_1,x_2},
\end{multline}
which is known as the \df{triple correlation} of \m{f}. 
Just as the bispectrum, 
provided that all Fourier matrices are invertible, 
the triple correlation also 
furnishes a complete set of left-invariant functionals of \m{f}. 
In fact, Kakarala derives the concept of bispectrum from the triple correlation, 
and not the other way round. 

\subsection{Computational considerations}

The price to pay for the completeness of the bipsectrum and the triple correlation is their inflated size. 
Letting \m{\absG} denote the cardinality of \m{G}, the triple correlation consists of 
\m{\absG^2} scalars. By the unitarity of the Fourier transform the total number of entries in the bispectrum matrices is the same. 
The symmetry \m{b(x_1,x_2)=b(x_2,x_1)} and its counterpart \m{\h b(\rho_1,\rho_2)=\h b(\rho_2,\rho_1)} 
(up to reordering of rows and columns) 
reduces the number of relevant components to \m{\absG(\absG+1)/2}, 
but to keep the analysis as simple as possible we ignore this constant factor. 
Another technical detail that we ignore for now is that if all we need is a 
complete set of invariants, then we can use the matrices  
\m{\h b\br{\rho_1,\rho_2} C_{\rho_1\rho_2}} instead of \m{\h b\br{\rho_1,\rho_2}},  
sparing us the cost of multiplying by \m{C_{\rho_1\rho_2}^\dag} 
in \rf{eq: bispectrum}. 

The cost of computing the bispectrum involves two factors: the cost of the Fourier transform, 
and the cost of multiplying the various matrices in \rf{eq: bispectrum}. 
Recent years have seen the emergence of fast Foruier transforms for a series of non-commutative groups, 
including the symmetric group \citep{Clausen89}, wreath product groups \citep{Rockmore95}, 
and others \citep{MaslenRockmore97}.  
Typically these algorithms reduce the complexity of Fourier transformation (and inverse Fourier transformation) from \m{O(\absG^2)} to 
\m{O(\absG \log^k\nts \absG)} scalar operations for some small integer \m{k}. 
In the following we assume that for whatever group we are working on, an \m{O(\absG \log^k\nts \absG)} transform is available. 
This makes matrix multiplication the dominant factor in the cost of computing the bispectrum. 

The complexity of the naive approach to multiplying two \m{D}-dimensional matrices is \m{D^3}. 
Using the well known fact that 
\m{\sum_{\rho\in\Rcal}d_\rho^2=\absG} (which can be regarded as a corollary to the Fourier transform 
being an invertible linear map\ignore{\m{L(G)\to \bigoplus_{\rho\in\Rcal} \CC^{d_\rho\times d_\rho}}}), 
the cost of computing \rf{eq: bispectrum} for given \m{(\rho_1,\rho_2)} is \m{O(d_{\rho_1}^3 d_{\rho_2}^3)}. 
Summing over all \m{\rho_1,\rho_2\in\Rcal} and 
assuming that \m{\sum_{\rho\in\Rcal} d_{\rho}^3} grows with some power \m{1< \theta <2} of \m{\absG} 
(note that even \m{\sum_{\rho\in\Rcal} d_{\rho}^4} grows with at most \m{\absG^2})  
gives an overall complexity bound of \m{O(\absG^{2\theta})}. 
Thus, the critical factor in computing the bispectrum is matrix multiplication and not 
the fast Fourier transform. 

As for the triple correlation, \rf{eq: triple correlation} involves an explicit sum over \m{G}, which is to be computed for 
all \m{\absG^2} possible values of \m{\br{x_1,x_2}}, giving a  
total time complexity is \m{O(\absG^3)}. 

One should note that when unique identifiability of \m{f} is not an absolute necessity, 
it is possible to truncate \m{\h f} according to some ``low pass filtering'' scheme 
in the interest of computational efficiency. 
This is easy to implement in the bispectrum, but not in the triple correlation. 

\section{The skew spectrum}\label{sec: skew}

The computational demands of the 
triple correlation 
and the bispectrum 
are of serious concern in applications. 
The former involves an explicit summation over \m{G} for each \m{(x_1,x_2)} pair, while 
the latter involves multiplying together large matrices and the non-trivial issue of 
computing Clebsch Gordan coefficients. In this section we derive a third, unitarily equivalent, set of invariants, 
which, in some sense are a combination of the two. 

Our starting point is the collection of functions  \m{r_z\colon G\to \CC} indexed by \m{z\in G} and defined  
\begin{equation*}
  r_z(x)=f(x)\, f(xz).
\end{equation*}
Introducing the concept of  \df{left \m{z}-diagonal slice} \m{g\Lslice{z}(x)=g(x,zx)}
and \df{right \m{z}-diagonal slice} \m{g\nts\Rslice{z}(x)=g(x,xz)} 
of a general function  \m{g\colon G\times G\to\CC}, and recalling 
our previous definition \m{u(x_1,x_2)=f(x_1)\ts f(x_2)}, 
we may write \m{r_z=u\Rslice{z}}~. 
The  \m{\h r_z(\rho)} Fourier components are \m{\rho}-covariant 
(with respect to the action of \m{G} on \m{f}), since  
\begin{equation*}
  \h r^{\,t}_{\nts z}(\rho)= 
  \sum_{x \in G} \rho(x)\, f(t^{-1}x)\, f(t^{-1}x z)
  =\sum_{x \in G} \rho(tx)\, f(x)\,f(xz)
  =\rho(t)\cdot \h r_{\nts z}(\rho). 
\end{equation*}
This immediately gives rise to the left-translation invariant matrices 
\begin{equation}\label{eq: skew spectrum}
  \h q_z(\rho)=\h r_z(\rho)^\dag\cdot \h f(\rho) \qqquad\rho\in\Rcal.  
\end{equation}
The collection of matrices \m{\cbr{\h q_z(\rho)}_{\rho\in\Rcal, z\in G}} we call the \df{skew spectrum} of \m{f}.
 
To see that the skew spectrum and the bispectrum are unitarily equivalent, 
it is sufficient to observe that by the convolution theorem 
\begin{multline*}
  q_z(x)=\sum_{y\in G} \mirr{r}_z(x y^{-1})^\ast\, f(y)
  =\sum_{y\in G} r_{\nts z}(y x^{-1})^\ast\, f(y)\\
  =\sum_{y\in G} f(y x^{-1})^\ast\,f(y x^{-1}\nts z)^\ast\, f(y)
  =b(x,z^{-1} x)
  =b\Lslice{z^{-1}}(x),  
\end{multline*}
and that \m{\cbr{b\Lslice{z}\,}_{z\in G}} together make up the triple correlation \m{b}. 
We thus have the following theorem. 

\begin{thm}\label{thm: skew complete}
  Let \m{f} and \m{f'} be complex valued functions on a finite group \m{G} 
  and let \m{\Rcal} be a complete set of inequivalent irreducible representations of \m{G}. 
  Assume that \m{\h f(\rho)} is invertible for all \m{\rho\in\Rcal}.  
  Let the skew spectrum \m{\h q_z} be defined as in \rf{eq: skew spectrum}. 
  Then \m{f'=f^t} for some \m{t\in G} if and only if \m{\h q'_z(\rho)=\h q_z(\rho)} for all \m{z\in G} and all \m{\rho\in\Rcal}. 
\end{thm}

From a computational perspective, the skew spectrum involves \m{\absG\<+1} separate Fourier transforms followed by 
\m{\absG} sequences of multiplying \m{\cbrN{\CC^{d_\rho\times d_\rho}}_{\rho\in\Rcal}} matrices. 
Using the notations of the previous section, the cost of the former is \m{O(\absG^2 \log^k\nts \absG)}, 
while the cost of the latter is \m{O(\absG^{\theta+1})}, improving on both the triple correlation and the bispectrum.  

\section{Homogeneous spaces}\label{sec: homo}

Real world problems often involve functions on homogeneous spaces of groups as opposed to functions on 
the groups themselves. 
Recall that a \df{homogeneous space} of a finite group \m{G} is a set \m{S} on which \m{G} acts by \m{s\mapsto x(s)}  
in such a manner that for any given \m{{s_0\in S}},~ \m{\setof{x(s_0)}{x \in G}} sweeps out the entire set. 
The group elements stabilizing \m{s_0} form a subgroup \m{H}, making it possible to identify  
any \m{s\tin S} with some coset \m{xH}. 
Thus, \m{S} itself is identified with the \df{quotient space} \m{G/H} of left cosets of \m{H} in \m{G}. 
A \df{transversal} for \m{G/H} is a set containing  
exactly one group element from each \m{xH} coset. 
By abuse of notation we use the symbol \m{G/H} for the transversals as well as the quotient space.
The right cosets \m{\cbr{Hx}_{x\in G}} also form a quotient space 
(with respect to the right action of \m{G}) 
and this we denote \m{H\tts\backslash\nts G}. 
We examine functions on such spaces in the next section. 

We denote the space of complex valued functions on \m{S=G/H} by \m{L(G/H)}. 
Any \m{f\in L(G/H)} extends naturally to a \df{right} \m{H}-\df{invariant} function \m{f\nts\tup{G}{}} on \m{G} by 
\m{  f\nts\tup{G}{}\!\!(x)=|H|^{-1} f(xH)}.
Conversely, any \m{g\in L(G)} 
may be restricted to \m{L(G/H)} by 
\m{  g\tdown{G/H}{}(xH)=\sum_{h\in H} g(xh)}.  
In accordance with this correspondence, for \m{z\in G} we define the \df{left-translate} of 
\m{f\in L(G/H)} as the function \m{f^z\in L(G/H)} given by 
\m{f^z(xH)=f(z^{-1}xH)}, and the Fourier transform of \m{f} as 
\begin{equation}\label{eq: FT homo}
  \h f(\rho)=\sum_{x\in G} \rho(x)\,f(xH)=\sqbbigg{\:\sum_{x\in G/H} \rho(x)\,f(xH)\:} \cdot  \sqbbigg{\:\sum_{h\in H} \rho(h)\:} .  
\end{equation}
Clearly, \m{\h f(\rho)=\h{f\tup{G}{}}(\rho)}, and \m{\h{f^z}(\rho)=\rho(z)\cdot \h{f}(\rho)} remains valid for functions on \m{G/H}.
\ignore{
with inverse transform
\begin{equation*}
  f(xH)=\ovr{\abs{G}} \sum_{\rho\in\Rcal} d_\rho\: \tr\brbig{\rho(x) \cdot \h f(\rho)} \qqquad g\in G. 
\end{equation*}
}

We denote by \m{\rho\tdown{H}{}} the restricted representations \m{\rho\colon H\to\CC^{d_\rho\times d_\rho}}  
given by \m{\rho\tdown{H}{}(h)=\rho(h)}. 
While \m{\rho} is always an irreducible representation of \m{G}, in general \m{\rho\tdown{H}{}} is reducible, i.e., 
\begin{equation}\label{eq: adapted}
  \rho\tdown{H}{}(h)=U_\rho\sqbbigg{\bigoplus_{\eta} \eta(h)} U_\rho^\dag \qqquad \qqquad h\in H,
\end{equation}
where \m{U_\rho} is a unitary matrix and 
\m{\eta} runs over some well-defined subset of irreducible representations of \m{H}, possibly with repeats. 
If \m{U_\rho=I} we say that \m{\rho} is \df{adapted} to \m{H} or, equivalently, that it is expressed in a 
\df{Gelfand-Tsetlin} basis with respect to \m{G/H}. 
In the following we assume that each \m{\rho\in\Rcal} is not only unitary but also \m{H}-adapted. 
It is possible to show that one may always choose \m{\Rcal} so as to satisfy these conditions. 

For the trivial representation \m{\eta_{\textrm{triv}}(h)=1} we have \m{\sum_h \eta_{\mathrm{triv}}(h)=\abs{H}}, hence   
by the unitarity of the Fourier transform over \m{H} for any other irreducible we must 
have \m{\sum_h \eta(h)=0}. 
Plugging \rf{eq: adapted} (with \m{U_\rho=I}) back in \rf{eq: FT homo} then shows that the Fourier 
transform of functions in \m{L(G/H)} have a very special form: only those columns of the \m{\h f(\rho)} 
matrices are non-zero which correspond to the trivial representation of \m{H} in \rf{eq: adapted}. 

The consequences of this type of sparsity are two-fold. 
On the one hand, except for the trivial case \m{H=e}, the sparsty implies that the non-singularity 
condition required for Kakarala's completeness result and derived results such as Theorem \ref{thm: skew complete} 
are \emph{always} violated when working over homogeneous spaces. 
The bispectrum, triple correlation and skew spectrum remain useful invariants, but we can no longer trust that 
they uniquely determine the original function up to translation. 
On the other hand, the sparsity suggests that the invariants can be computed much faster than when we were working 
over the full group \m{G}. 

\ignore{
Except for the trivial case \m{H\<=e},  
the Fourier transform \m{\mathfrak{F}\colon f\mapsto \h f} is no longer unitary, since it is no longer surjective. 
Furthermore, if for some \m{\rho\in\Rcal} we have  
\m{  \sum_{h\in H} \rho(h)=0},
the corresponding component of \m{\h f} will vanish. 
We denote by \m{\Rcal_{G/H}} the subset of \m{\Rcal} for which this 
does \emph{not} occur. 
In many cases of interest \m{\Rcal_{G/H}} is a fairly small set. 

Another simplifying factor when working over homogeneous spaces is that 
even for \m{\rho\in \Rcal_{G/H}} the Fourier matrices might have vanishing columns. 
We say that a representations \m{\rho} of \m{G} is \df{adapted} to the subgroup \m{H} 
(equivalently, that it is expressed in a \df{Gelfand-Tsetlin basis}), if on restriction to \m{H} 
it splits in the form 
\begin{equation}\label{eq: adapted}
  \rho(h)=\bigoplus_{\eta} \eta(h) \qqquad \qqquad h\in H,
\end{equation}
where \m{\eta} runs over some subset of irreducible representations of \m{H}, possibly with repeats. 
From the unitarity of the Fourier transform it is immediate that except for the trivial representation 
\m{\eta_{\text{triv}}(h)\equiv 1}, we have \m{\sum_{h\in H}\eta(h)=0}. 
Consequently, 
only those columns of \m{\h f(\rho)} will have non-vanishing entries 
which correspond to trivial \m{\eta} in \rf{eq: adapted},  
effectively reducing the Fourier matrices to a collection of column vectors. 
}


Letting \m{n_\rho} be the multiplicity of the trivial representation in the summation on the right hand side of  
\rf{eq: adapted}, each bispectrum component (with the final \m{C_{\rho_1\rho_2}^\dag} omitted)
\begin{equation}\label{eq: bispectrum homo}
  \h b(\rho_1,\rho_2)=\brbig{\h f(\rho_1) \otimes \h f(\rho_2)}^\dag\:
  C_{\rho_1\rho_2}\:
  \bigoplus_{\rho\in\Rcal} {\h f(\rho)}^{\oplus m_\rho} 
\end{equation}
will only have \m{n_{\rho_1} n_{\rho_2}} non-zero rows and a similary small number of non-zero columns.


Since \m{y(hx)^{-1}\nts H=yx^{-1}hH=yx^{-1}H} for any \m{h\in H}, the triple correlation 
\begin{equation*}
  b(x_1,x_2)=\sum_{y\in G} f(y x_1^{-1}H)^\ast\,f(y x_2^{-1}H)^\ast\,f(yH)
\end{equation*}
will be a left \m{H\times H}-invariant function, i.e., effectively \m{b\in L(H\backslash G\times H\backslash G)}. 
The summation still extends over all of \m{y\in G}, however, giving a total complexity of \m{\absG^3/\abs{H}^2}.



As for the skew spectrum, the first fact to note is that the \m{z} index can be restricted to 
one element from each \m{HxH=\setof{h_1 x h_2}{h_1,h_2\in H}} \df{double coset}. 
A transversal for such double cosets we denote \m{H\backslash G/H}. 
\begin{thm}\label{thm: skew restrict}
  Let \m{H} be  subgroup of a finite group \m{G} and let \m{f\in L(G/H)}. 
  Then the skew spectrum \m{\h q_f}  is uniquely determined 
  by its subset of components \m{\setof{\h q_z(\rho)}{\rho\in\Rcal,~z\in H\backslash G/H}}. 
\end{thm}
\begin{pf}
  For any \m{h\in H}, 
  \begin{equation*}
    r_{zh}(x)=f(xzhH)f(xhH)=f(xzH)f(xH)=r_z(x), 
  \end{equation*} 
  so \m{\h q_{zh}(\rho)=\h r_z(\rho)^\dag \h f(\rho)=\h r(\rho)^\dag \h f(\rho)=q_z(\rho)}. 
  Now for \m{h'\in H}
  \begin{equation*}
    r_{h'z}(x)=f(xh'zH)f(xH)=f(xh'zH)f(xh'H)=r_{z}^{({h'}^{-1})}(x),  
  \end{equation*}
  where \m{r_{z}^{(t)}} denotes the right-translate  \m{r_{z}^{(t)}(x)=r_z(xt^{-1})}.  
  Thus, by the right-translation property 
  \begin{equation*}
  \h {f^{(t)}}(\rho)=\sum_{x \in G} \rho(x)\, f(xt^{-1})
  =\sum_{x \in G} \rho(xt)\, f(x)
  =\sqbbigg{\sum_{x \in G} \rho(x)\, f(x)}\rho(t)
  =\h f(\rho) \cdot \rho(t) 
  \end{equation*}
  of the Fourier transform,
  \begin{equation*} 
    {\h q_{h'z}(\rho)=(\h r_z(\rho)\, \rho({h'}^{-1}))^\dag \h f(\rho)=
    \rho(h')\: \h r_z(\rho)^\dag\ts \h f(\rho) = \rho(h')\: \h q_{z}(\rho)}, 
  \end{equation*} 
  so \m{\h q_{h'z}(\rho)} and \m{\h q_z(\rho)} albeit not equal, are related by an invertible linear mapping. 
  \hfill 
\end{pf} 
\\ \\
What the best way to compute \m{\h q_z(\rho)} is strongly depends on \m{G} and \m{H}. 
One possible algorithm is based on the decomposition 
\begin{equation*}
  \h f(\rho)=\sum_{y\in G/H} \rho(y) \bigoplus_\eta \sum_{h\in H} \eta(h)\, f(yh)
  =\sum_{y\in G/H} \rho(y) \bigoplus_\eta \h f_y(\eta) 
\end{equation*}
where the \m{f_y} are functions on \m{H} defined \m{f_{y}(h)=f(yh)} and 
\m{\h f_y} are their Fourier transforms
\begin{equation*}
  \h f_y(\eta)=\sum_{h\in H} \eta(h) f_y(h). 
\end{equation*}
This is akin to a Clausen-type ``matrix separation of variables'' fast Fourier transform (FFT) tailored to the 
\m{G>H>1} chain. 
Because \m{f\in L(G/H)}, all but the \m{\h f_y(\eta_{\mathrm{triv}})} component of \m{\h f_y} vanish, and  
\m{\h f} can be computed efficiently by restricting the FFT 
to the descendants of the \m{\h f_y(\eta_{\mathrm{triv}})} components  
(see the section on ``partial Fourier transforms'' in the documentation to \texttt{SnOB} \citep{Snob} for details). 
Letting \m{g} be the right-translated function 
\m{g(x)=f^{(z^{-1})}(x)=f(xz)} we can then compute \m{\h g (\rho)=\h f(\rho)\, \rho(z^{-1})} for each \m{\rho} and 
run the FFT backwards to get \m{\h g_y(\eta)} for each \m{y\in G/H}. Multiplying each of these by the corresponding 
\m{f(y)} yields \m{\h g_y(\eta)f(y)}, which is exactly the Fourier transform of \m{r_z} restricted to \m{yH}, so 
one more forward transform yields \m{\h r_z(\rho)}, which we can plug directly into \rf{eq: skew spectrum}. 
The exact complexity of this procedure depends on the groups \m{G} and \m{H}.

\ignore{
On \m{G/H} the skew spectrum \m{\h q_z(\rho)=\h r_z(\rho)^\dag\cdot \h f(\rho)} 
restricts to \m{\Rcal\downarr_{G/H}}, and in 
\begin{equation*}
  r_z(x)= f(xH)\,f(xzH)
\end{equation*}
it is apparent that the index \m{z} may also be restricted to \m{G/H}. 
The function \m{r_z} is still a function on the entire group \m{G}, however, 
and thus computing \m{\h r_z} requries a Fourier transform over \m{G}.  
To address this, we first note that since \m{\h f(\rho)} is restricted to the subspace of 
Fourier transforms of functions from \m{L(G/H)}, \m{r_z(\rho)} may be projected to 
the same subspace without changing \m{\h q(\rho)}. 
In other words, 
\begin{equation*}
  \h q_z(\rho)=\h r_z(\rho)^\dag\cdot \h f(\rho)=\h r^{\circ}_z(\rho)^\dag\cdot \h f(\rho),
\end{equation*}
where \m{\h r^{\circ}_z} is the Fourier transform of the \m{L(G/H)} function 
\begin{equation*}
  r_z\tdown{G/H}{}(xH)=\sum_{h\in H} f(xhH)\,f(xhzH)=
  f(xH) \sum_{h\in H} f(xhzH)=f(xH)\,f_z\tdown{G/H}{}(xH),
\end{equation*}
where \m{f_z(x)=f(xzH)}. It remains to compute \m{f_z\tdown{G/H}{}}, which we do 
by additional Fourier transforms over \m{G/H}.  

First note that \m{f_e=f\tup{G}{}\;\in L(G/H)}, so \m{\h f_e} is easy to compute. 
Secondly, note that \m{f_z} is the \df{right-translate} of \m{f_e} by \m{z^{-1}}, i.e., 
\m{f_z(x)=f_e(xz)}. Now in exact analogy with \rf{eq: translation}, 
\begin{equation*}
  \h f_z(\rho)=\sum_{x\in G} \rho(x)\,f_e(xz)=\sum_{x\in G} \rho(xz^{-1})\,f_e(x)=
  \sqbbigg{\:\sum_{x\in G} \rho(x)\,f_e(x)\:}\rho(z^{-1})=\h f_e(\rho)\, \rho(z)^\dag, 
\end{equation*}
from which \m{f_z} may be computed by inverse Fourier transformation over \m{G/H}. 

In summary, the computation leading to \m{\h q_z} may be written as 
\begin{equation*}
  \h q_z=\Fcal\brbig{\:f \cdot \Fcal^{-1}\sqbbig{\:\Fcal(f)\odot\br{\rho(z)^\dag}_\rho\: }\:} \odot \Fcal(f),
\end{equation*}
where \m{\Fcal} is the Fourier transform over \m{G/H}, \m{\Fcal^{-1}} is its inverse, 
\m{\odot} denotes taking the product of the corresponding Fourier matrices, and
\m{\cdot} denotes the (pointwise) product of functions in \m{G/H}. 
}

\ignore{
and 
\begin{equation*}
  f_z(xH)=\sum_{h\in H} f(xhzH) 
\end{equation*}
is itself a function in \m{L(G/H)}. 
It remains to compute \m{f_z} without explicitly executing the summations over \m{h}. 
This is done by noting that 
\begin{equation*}
  f_z\tup{G}{}\!(x)=\sum_{h\in H} f(xhzH) 
\end{equation*}
}
\ignore{
\begin{equation*}
  \h r'_z(\rho)=\sum_{x\in G/H} \rho(x)\, f(xH) \sum_{h\in H} f(xhzH)
\end{equation*}

\begin{equation*}
  \br{f\tup{G}{}}^{z^{-1}}\tdown{G/H}{}(x)
\end{equation*}
}
\ignore{
Similarly to \rf{eq: defh}, we let 
\begin{equation}\label{eq: homo g}
  g_z(x)=f(xH)\,f(xzH), 
\end{equation}
but note that for any \m{h\in H},~ \m{g_{zh}=g_z}, hence it is sufficient to define \rf{eq: homo g} 
for a complete set \m{Z=\cbr{\seq{z}{\abs{\!G/H\!}}}} of coset representatives of \m{G/H}. 
Now there are two distinct cases to consider. 
If \m{H} is a normal subgroup of \m{G}, then for any \m{h\in H},~\m{hz\in zH}, and hence  
\begin{equation*}
  g_z(xh)=f(xH)\,f(xhzH)=f(xH)\,f(xzH)=g_z(x). 
\end{equation*}
In this case \m{g_z} may itself be regarded as a function on \m{G/H}, which can potentially significantly 
simplify the computations. 
If \m{H} is not a normal subgroup, then 
\begin{equation*}
  g_z(xh)=g_{hz}(x).
\end{equation*}
\begin{equation*}
  \h g_z(\rho)=\sum_{x\in G} \rho(x) f(xH)\,f(xzH)
\end{equation*}
\m{x'\in xH \cap z^{-1} x z H}
\begin{equation*}
  \h q_z(\rho)=
\end{equation*}
}

\clearpage
\section{Right invariance}

The \df{right-translate} of \m{f\colon G\to\CC} by \m{z\in G} is defined \m{f^{(z)}(x)=f(xz^{-1})}, 
and a functional \m{p\colon f\mapsto p(f)} is said to be \df{right-invariant} if \m{p(f^{(t)})=p(f)} 
for any \m{f} and any \m{t\in G}. 
Everything that we described in sections \ref{sec: bispectrum} and \ref{sec: skew} have natural right-invariant analogs.
The Fourier transform obeys \m{\h f^{(t)}(\rho)=\h f(\rho)\, \rho(t)} and 
the corresponding invariant power spectrum is \m{\h q^{R}(\rho)=\h f(\rho)\, \h f(\rho)^\dag}, 
which is the Fourier transform of the ``right-autocorrelation'' \m{q^{R}(x)=\sum_{y\in G}f(y)\,f(xy)}. 
Similarly, the ``right-bispectrum'' is 
\begin{equation*}
  \h b^R(\rho_1,\rho_2)=\brbig{\h f(\rho_1) \otimes \h f(\rho_2)}\:
  C_{\rho_1\rho_2}\:
  \sqbbigg{\:\bigoplus_{\rho\in\Rcal} {\h f(\rho)^\dag}^{\oplus m_\rho} \:}\:  C_{\rho_1\rho_2}^\dag, 
\end{equation*}
which is the Fourier transform of \m{b^R(x_1,x_2)=\sum_{y\in G}f(x_1y)\,f(x_2y)\,f(y)}.  
For the skew spectrum we must change \m{r_z} to \m{r^R_z(x)=u\Lslice{z}(x)=f(x)\,f(zx)}, 
and let \m{\h q^R_z(\rho)=\h r^{R}_z(\rho)\,\h f(\rho)^\dag}, which is the Fourier transform of 
\m{q^R_z=\sum_{y\in G} f(xy)\,f(zxy)\,f(y)}. 
The rest of the analysis goes through exactly as in the left-invariant case. 

Generalizing section \ref{sec: homo} is more interesting, 
because our choice of \m{G/H} over \m{H\tts\backslash\nts G} breaks the left-right symmetry.  
For notational convenience, instead of constructing right-invariant functionals of \m{f\in L(G/H)}, 
we discuss the analogous problem of left-invariant functionals of \m{f\in L(H\backslash G)}.  

Left-translating \m{f\in L(H\backslash G)} by \m{t\in H} leaves it invariant, 
effectively reducing the set of transformations of interest from \m{G} to \m{G/H}. 
In this sense we are in a easier situation than in section \ref{sec: homo}. 
Indeed, despite that the Fourier matrices of \m{f} are just as rank deficient as 
in the previous section, in this case Kakarala's completeness result is salvagable. 

\begin{thm}\cite[Theorem 3.3.6]{KakaralaPhD}
  Let \m{H} be any closed subgroup of a compact group \m{G}, and let \m{f\in L_1(H\backslash G)} be such that 
  for all \m{\rho\in\Rcal}, the matrix rank of \m{\h f(\rho)} is equal to the multiplicity of \m{\eta_{\textrm{triv}}} 
  in the decomposition of \m{\rho} into irreducible representations of \m{H}. Then \m{\h b_f=\h b_{f'}} for some 
  \m{f'\in L_1(G)} if and only if there exists some \m{t\in G} such that \m{f'=(f\tup{G}{}\,)^t}. 
\end{thm}
 
As a corollary, the skew spectrum will all so be complete. However, the skew spectrum is a much larger 
object than before because 
(a) there is no obvious way to restrict \m{z} like in Theorem \ref{thm: skew restrict} 
(b) instead of being columns-sparse, the Fourier matrices are row-sparse, hence as long as the decomposition 
of \m{\rho\tdown{H}{}} contains \m{\eta_{\mathrm{triv}}} with multiplicity at least one, \m{\h q_z(\rho)} will 
be a full \m{d_\rho\times d_\rho} matrix.

\section{Conclusions}

We have presented...

\section*{Acknowledgements}

I would like to thank Ramakrishna Kakarala for providing me with a hard copy of his thesis. 

\bibliography{bispectrumFFT}

\begin{thebibliography}{6}
\providecommand{\natexlab}[1]{#1}
\providecommand{\url}[1]{\texttt{#1}}
\expandafter\ifx\csname urlstyle\endcsname\relax
  \providecommand{\doi}[1]{doi: #1}\else
  \providecommand{\doi}{doi: \begingroup \urlstyle{rm}\Url}\fi

\bibitem[Clausen(1989)]{Clausen89}
M.~Clausen.
\newblock Fast generalized {F}ourier transforms.
\newblock \emph{Theor. Comput. Sci.}, pages 55--63, 1989.

\bibitem[Kakarala(1992)]{KakaralaPhD}
R~Kakarala.
\newblock \emph{Triple corelation on groups}.
\newblock PhD thesis, Department of Mathematics, UC Irvine, 1992.

\bibitem[Kakarala(1993)]{KakaralaTriple}
R.~Kakarala.
\newblock A group theoretic approach to the triple correlation.
\newblock In \emph{IEEE Workshop on higher order statistics}, pages 28--32,
  1993.

\bibitem[Kondor(2006)]{Snob}
Risi Kondor.
\newblock $\mathbb{S}_n$\texttt{ob}: a {C++} library for fast {F}ourier
  transforms on the symmetric group, 2006.
\newblock {A}vailable at \texttt{http://www.cs.columbia.edu/\~{ }risi/Snob/}.

\bibitem[Maslen and Rockmore(1997)]{MaslenRockmore97}
D.~Maslen and D.~Rockmore.
\newblock Generalized {FFT}s --- a survey of some recent results.
\newblock In \emph{Groups and Computation {II}}, volume~28 of \emph{{DIMACS}
  Ser. Discrete Math. Theor. Comput. Sci.}, pages 183--287. AMS, Providence,
  RI, 1997.

\bibitem[Rockmore(1995)]{Rockmore95}
D.~Rockmore.
\newblock Fast fourier transforms for wreath products.
\newblock \emph{J. Applied and Computational Harmonic Analysis}, 2:\penalty0
  279--292, 1995.

\end{thebibliography}
\bibliographystyle{plainnat}
\end{document}